\documentclass[twoside,a4paper,reqno,11pt]{amsart}
\usepackage{amsfonts, amsbsy, amsmath, amssymb, latexsym}
\usepackage{mathrsfs,array}
\usepackage[top=30mm,right=30mm,bottom=30mm,left=30mm]{geometry}
\usepackage{stmaryrd}
\usepackage{bm}

\usepackage{hyperref}

\def\d{\delta}

\def\G{\Gamma}
\newcommand{\F}{\mathbb{F}}

\newcommand{\N}{\mathbb{N}}
\newcommand{\e}{\epsilon}
\renewcommand{\bf}{\textbf} 
 \renewcommand{\O}{\Omega}

 \renewcommand{\to}{\rightarrow}

\newcommand{\la}{\langle}
\newcommand{\ra}{\rangle}

\makeatletter
\newcommand{\imod}[1]{\allowbreak\mkern4mu({\operator@font mod}\,\,#1)}
\makeatother

\newtheorem{theorem}{Theorem}
\newtheorem*{conj*}{Conjecture}

\newtheorem{propn}[theorem]{Proposition}

\newtheorem{corol}[theorem]{Corollary}

\theoremstyle{definition}

\begin{document}

 \author{Martin W. Liebeck}
\address{M.W. Liebeck, Department of Mathematics,
    Imperial College, London SW7 2BZ, UK}
\email{m.liebeck@imperial.ac.uk}

\author{Aner Shalev}
\address{A. Shalev, Institute of Mathematics, Hebrew University, Jerusalem 91904, Israel}
\email{shalev@math.huji.ac.il}

\title{Girth, words and diameter}

\begin{abstract}
We study the girth of Cayley graphs of finite classical groups $G$ on random sets of generators. Our main tool is an essentially best possible  bound we obtain on the probability that a given word $w$ takes the value 1 when evaluated in $G$ in terms of the length of $w$, which has
additional applications. We also study the girth of random directed Cayley graphs of symmetric groups, and the relation between the girth and the diameter of random Cayley graphs of finite simple groups.
 \end{abstract}


\subjclass[2010]{Primary 20D06; Secondary 20P05, 05C80, 05C12}

\thanks{We are grateful to Sean Eberhard for some enlightening comments.
The second author acknowledges the support of ISF grant 686/17, BSF grant 2016072 and the Vinik chair of mathematics which he holds.}

\maketitle

\section{Introduction}

The girth of a graph (resp. directed graph) is the minimal length of a cycle (resp. directed cycle) in the graph. The girth of finite $k$-regular graphs has been studied extensively, with a particular focus on graphs of large girth -- see for example \cite {ES}, \cite{Mar}. Note the trivial upper bound of $2\log_{k-1}{v}+1$ for the girth, where $v$ is the number of vertices and $k>2$.

 In \cite{GHSSV} the girth of random Cayley graphs of various families of groups was studied, and large girth results were established.
For a finite group $G$ and a sequence $S$ of elements $g_1,\ldots ,g_k$ of $G$, let $\G(G,S)$ (resp. $\G^*(G,S)$) denote the associated undirected (resp. directed) Cayley graph.  Corollary 2 of \cite{GHSSV} asserts that for finite simple groups $G$, the girth of $\G(G,S)$ for $k$ random generators tends to $\infty$ almost surely as $|G| \to \infty$. Also \cite[Thm. 4]{GHSSV} shows that for groups $G$ of Lie type of bounded rank, the girth is $\O(\log |G|)$, while \cite[Thm. 3]{GHSSV} asserts that for $G = S_n$, the girth is at least $\O((\log |G|)^{1/2})$ almost surely. An error in the proof of the latter result was recently pointed out in \cite{Eber}, and a slightly weaker bound of the form  $\O(n^{1/3})$ was obtained.

The random girth of classical groups $G$ of unbounded rank has apparently remained unexplored.
Denote by $Cl_n(q)$ a simple classical group over $\F_q$ with natural module of dimension $n$.
 Our first result shows in particular that if the underlying field $\F_q$ has bounded size, then the random girth of such groups is at least
$\O( (\log |G|)^{1/2})$.

\begin{theorem}\label{girth} There exists an absolute constant $N$, and for each integer $k\ge 2$ and prime power $q$ a positive real number $b = b(q,k)$ with the following property. Let $G = Cl_n(q)$ with $n \ge N$, and let $S$ be a sequence of $k$ independently chosen random elements of $G$.
Then as $|G|\to \infty$, the girth of the Cayley graph $\G(G,S)$ exceeds $b\sqrt{\log |G|}$ almost surely.
\end{theorem}

In particular, for bounded $q$ and $k$ the girth is almost surely $\O(\sqrt{\log |G|})$.
Also the proof shows that the girth exceeds $Bn$ for some positive constant $B = B(k)$; in fact we can take
$B=\frac{1}{7(1+2 \log_2{(2k-1)})}$.

Our next result concerns the girth of directed Cayley graphs of symmetric groups (namely, the minimal length
of a directed cycle in the graph).

\begin{theorem}\label{pos}
 Fix an integer $k\ge 2$, and let  $S$ be a sequence of $k$ independently chosen random elements of $G=S_n$. Then the girth of  $\G^*(G,S)$ almost surely exceeds $c\sqrt{\log_k |G|}$, where $c$ is a positive absolute constant.
\end{theorem}

Thus the girth of a random directed Cayley graph of $S_n$ is at least $\O(\sqrt{n \log n})$.

The diameter of Cayley graphs of finite simple groups (with explicit or with random generators) has also attracted considerable attention --
see for instance \cite{BS}, \cite{BGT}, \cite{PS}, \cite{HS}, \cite{HSZ}, \cite{BT}, \cite{BY} and the references therein. Clearly if $d$ and $g$ are the diameter and girth respectively, then a trivial lower bound for $d$ is $\lfloor \frac{g}{2}\rfloor$, and there is interest in finding families of graphs for which $d$ is bounded in terms of $g$. According to \cite{AB}, a family of graphs is $dg$-{\it bounded} if the ratio $\frac{d}{g}$ is bounded. The focus is on graphs of {\it large} girth, meaning that the girth is $\O(\log |G|)$, where $G$ ranges over the ambient family of groups. The main result of \cite{AB} is a construction of certain Cayley graphs of $SL_n(p)$ with respect to two explicit generators, where $n$
is fixed, which are of large girth and $dg$-bounded.

It follows from part (i) of the next result that such families of Cayley graphs exist for all groups of Lie type of bounded rank.
Parts (ii) and (iii) bound the diameter in terms of (nonlinear) functions of the girth almost surely for other families of finite simple
groups.

\begin{propn}\label{DGbd} Fix $k \in \N$.
Let $G$ be a finite simple group and let $S$ be a sequence of $k$ independently chosen random elements of $G$.
Let $d, g$ be the diameter and girth of $\G(G,S)$ respectively. Suppose $|G| \to \infty$.
\begin{itemize}
\item[(i)] If $G$ is of Lie type of bounded rank, then $\G(G,S)$ has large girth and is $dg$-bounded (i.e. $d=O(g)$) almost surely.
\item[(ii)] If $G = Alt_n$, then almost surely $d \le g^6 (\log g)^c$ for some absolute constant $c$.
\item[(iii)] If $G = Cl_n(q)$ is a classical group of dimension $n$ with $q$ bounded, then
$d \le C^{g (\log g)^3}$ almost surely, for some constant $C=C(k) > 1$.
\end{itemize}
\end{propn}

The proof of part (i) is rather short, modulo the deep results in \cite{BGT, PS} (which are also used in \cite{AB}).
Parts (ii) and (iii) require results from \cite{HSZ} and  \cite{BY} respectively.

The bound in part (iii) above seems far from best possible, and it would be nice to obtain a polynomial bound in this case too.
Such a bound would follow from Theorem \ref{girth}, together with  Babai's conjecture (so far unproved) that the diameter of any connected Cayley graph of a (nonabelian) finite simple group $G$ is at most $(\log{|G|})^c$,
for some absolute constant $c$ (see \cite[1.7]{BS}). In fact, a polynomial bound in part (iii) of Proposition \ref{DGbd} would already follow if the latter bound
on the diameter holds almost surely for random Cayley graphs of classical groups.

The proofs of Theorems \ref{girth} and \ref{pos} rely on the study of the probability $P_G(w)$ that a word $w=w(x_1,\ldots,x_k)$ in the free group $F_k$ takes the value 1 when we substitute a sequence $S$ of $k$ independently chosen random elements of $G$ for $x_1,\ldots,x_k$. It is an elementary observation that
\[
{\bf P} (\hbox{girth}(\G(G,S)) \le L) \le  \sum_{|w|\le L} P_G(w)
\]
where $|w|$ denotes the length of $w$ (see \cite[Sec. 2]{GHSSV}). The study of $P_G(w)$ is an important part of the theory of word maps, with a particular focus on finite  simple groups $G$. In \cite[Thm. 3]{DPSS} it is shown that if $w\ne 1$ then $P_G(w) \to 0$ as $|G|\to \infty$; and \cite[Thm. 1.1]{LS} shows that for every $w\ne 1$ there exist $\d = \d(w)>0$ and $N = N(w)$ such that  $P_G(w) \le |G|^{-\d}$ provided $|G|>N$. The next result gives an explicit and close to best possible value for the constant $\d(w)$ in terms of the length of $w$.

\newpage

\begin{theorem}\label{prob}
For any $\e > 0$, there exists $c =c(\e)>0$ with the following property.
Let $\ell \in \N$ and let $G = Cl_n(q)$ be a classical group of dimension $n \ge c\ell$. Then, for any reduced word  $w \in F_k$
of length $\ell$, we have
\[
P_G(w) \le |G|^{-\frac{1}{(2+\e)\ell}}.
\]
\end{theorem}

In fact the proof gives $c(\e) =  4\,(1+\frac{2}{\e})$ for $G = SL_n(q)$, and
$c(\e) =  7\,(1+\frac{2}{\e})$ for the other classical groups. For more detailed bounds on $P_G(w)$ see Section 2.
Theorem \ref{prob} improves a bound of the form $P_G(w) \le  |G|^{-\frac{1}{1800\ell^2}+o(1)}$ obtained in the proof of \cite[Thm. 1.1]{LS}.

As claimed above, Theorem \ref{prob} is close to being best possible; indeed, \cite[Thm. 1.4]{LiSh3} shows that, for a fixed power word $w = x_1^{\ell}$ and
$G = Cl_n(q)$ with $n \to \infty$ we have $P_{w}(G) = |G|^{-\frac{1}{\ell} + o(1)}$. It seems an interesting and challenging problem to improve the upper bound in Theorem \ref{prob} to $P_G(w) \le |G|^{-\frac{1}{(1+\e)\ell}}$.

The key to the proof of Theorem \ref{pos} is the following result, which is of some independent interest.
Recall that a word $w \in F_k$ is said to be {\it positive} (or {\it a semigroup word}) if it does not involve inverses
of the generators of $F_k$.

\begin{propn}\label{posprop}
Let $w \in F_k$ be a positive word of length $\ell$. Then for all $n\in \N$, we have
\[
P_{S_n}(w) \le \left(\frac{2\ell}{n}\right)^{\lfloor \frac{n}{2\ell} \rfloor} \le (n!)^{-\frac{1}{2\ell} + o_n(1)}.
\]
\end{propn}

In fact the inaccurate proof of the above bound for $P_{S_n}(w)$ on p.106 of \cite{GHSSV}  becomes  accurate when $w$ is assumed to be a positive word. Proposition \ref{posprop} is essentially best possible; indeed for $w=x_1^{\ell}$ we have $P_{S_n}(w) = (n!)^{-\frac{1}{\ell} + o_n(1)}$
(see for instance \cite[2.17]{LiSh3}).

We conclude the introduction with applications of the two results above to representation varieties and subgroup growth
(cf. \cite[1.3, 1.4]{LS}).

\begin{corol}\label{rep} Let $\Gamma$ be a non-free group with $k$ generators, and let $\ell$ be the minimal length of a non-trivial relation
(in these generators) which holds in $\Gamma$. Then for every $\e > 0$ there exists $N = N(\ell, \e)$ such that the following hold for all $n \ge N$ and for any algebraically closed field $F$:
\begin{itemize}
\item[(i)] $\dim {\rm Hom}(\Gamma, GL_n(F)) \le (k- \frac{1}{(2+\e)\ell})n^2$;
\item[(ii)] $\dim {\rm Hom}(\Gamma, G) \le (k- \frac{1}{(2+\e)\ell})\dim G$, where $G$ is a simple algebraic group over $F$ of dimension $n$.
\end{itemize}
\end{corol}

Recall that $a_n(\Gamma)$ denotes the number of index $n$ subgroups of $\Gamma$.

\begin{corol}\label{growth} Let $\Gamma$ be a group with $k$ generators which satisfy some non-trivial positive relation.
Let $\ell$ be the minimal length of such a relation. Then $a_n(\Gamma) \le (n!)^{k-1-\frac{1}{2\ell} + o_n(1)}$.
\end{corol}

\section{Proof of Theorem \ref{prob}}

First we give the proof in the case $G= GL_n(q)$. Our method is inspired by Eberhard's proof in \cite{Eber}.
 Denote by $V = \F_q^n$ the underlying vector space. Assume $2\ell < n$.

Let $a_1,\ldots,a_k$ be free generators for $F_k$, and let $w = w_{\ell}\cdots w_1$, where each $w_i \in \{a_1^{\pm1},\ldots ,a_k^{\pm1}\}$. Let $g_1\ldots,g_k$ be a random sequence of elements of $G$.

Fix $v_1 \in V\setminus 0$, and define $v_1^0,\ldots v_1^{\ell}$ by
\[
\begin{array}{l}
v_1^0 = v_1, \\
v_1^j = w_j(g_1,\ldots ,g_k) (v_1^{j-1})\;\;(1\le j\le \ell).
\end{array}
\]
Call the sequence $v_1^0,\ldots v_1^{\ell}$ the {\it trajectory} of $v_1$.
Assume $v_1^0,\ldots ,v_1^{j-1}$ are linearly independent. Then
\begin{equation}\label{rest}
v_1^j \not \in {\rm Sp}\left( \{v_1^i\,(i\le j-1)\, |\, w_i=w_j \hbox{ or }w_{i+1}=w_j^{-1}\}\right),
\end{equation}
which excludes at most $q^{j-1}$ possibilities for $v_1^j$; all other vectors are equally likely as possibilities for $v_1^j$, since $w_j(g_1,\ldots ,g_k)$ is a random element.  Hence the conditional probability
\[
{\bf P}\left(v_1^j \in {\rm Sp}(v_1^0,\ldots ,v_1^{j-1})\,|\,v_1^0,\ldots ,v_1^{j-1}\right) \le \frac{q^j}{q^n-q^{j-1}}.
\]
It follows that
\[
{\bf P}\left(v_1^0,\ldots ,v_1^\ell \hbox{ lin. dep.}\right) \le \sum_{j=1}^\ell \frac{q^j}{q^n-q^{j-1}} \le
\frac{q+q^2+ \ldots + q^{\ell}}{q^n-q^{\ell-1}} < \frac{q}{q-1} \cdot \frac{q^{\ell}}{q^n-q^{\ell-1}}.
\]

Now suppose $v_1^0,\ldots ,v_1^\ell$ are given, and set $V_1 = {\rm Sp}(v_1^0,\ldots v_1^{\ell})$.
Pick $v_2 \not \in V_1$, and define the trajectory of $v_2$ to be
$v_2^0,\ldots v_2^{\ell}$ as above. Assuming $v_2^0,\ldots v_2^{j-1}$ to be linearly independent and also have span intersecting $V_1$ trivially, we have
\[
{\bf P}\left(v_2^j \in {\rm Sp}(v_2^0,\ldots v_2^{j-1} \cup V_1)\right) \le \frac{q^{\ell+j}}{q^n-q^{\ell+j-1}},
\]
and hence, arguing as above, we obtain
\[
{\bf P}\left(v_2^0, \ldots , v_2^{\ell}\hbox{ lin. dep.}\,|\, v_1^0,\ldots ,v_1^\ell  \right) \le
\sum_{j=1}^\ell  \frac{q^{\ell+j}}{q^n-q^{\ell+j-1}} < \frac{q}{q-1} \cdot \frac{q^{2\ell}}{q^n-q^{2\ell-1}}.
\]
Repeating this argument $m$ times, where $m \le \frac{n}{\ell}$  (choosing each $v_i$ ($1\le i\le m$) not in the span of the previous trajectories), we obtain
\[
{\bf P}\left(v_m^0, \ldots , v_m^{\ell} \hbox{ lin. dep.}\,|\, v_i^0,\ldots ,v_i^\ell \hbox{ for } i<m \right) <
 \frac{q}{q-1} \cdot \frac{q^{m\ell}}{q^n-q^{m\ell-1}}.
\]
If $w(g_1,\ldots,g_k)=1$, then $v_i^\ell = v_i$ for all $i$, and hence
\[
\begin{array}{ll}
P_G(w) & \le \prod_{i=1}^m {\bf P}\left(v_i^\ell = v_i\,|\,v_j^\ell=v_j \hbox{ for } 1\le j\le i-1\right) \\
 & < \prod_{i=1}^m  \frac{q}{q-1} \frac{q^{i\ell}}{q^n-q^{i\ell-1}} = (\tfrac{q}{q-1})^m\prod_{i=1}^m \frac{1}{q^{n-i\ell}-q^{-1}}.
\end{array}
\]
Set $m=\lfloor\frac{n}{\ell} \rfloor$. Define
\[
a(q) = \prod_{i=1}^m \frac{q^{n-i\ell}}{q^{n-i\ell}-q^{-1}} = \prod_{i=1}^m \frac{1}{1 - q^{-(n-i\ell) -1}}.
\]
We may and shall assume $\ell \ge 2$, since for $w$ of length $1$ we have $P_G(w) = |G|^{-1}$ for all finite
groups $G$. Clearly $a(q) \le a(2) \le a$, where
\[
a := \prod_{i = 0}^{\infty} \frac{1}{1 - 2^{-(2i+1)}} = \prod_{i = 0}^{\infty} (1 + \frac{1}{2^{2i+1}-1})
< 2.3749,
\]
where the last inequality is easily verified by computing the sum for $i \le 6$ and bounding its tail.
It follows that
\[
P_G(w) \le a \cdot (\tfrac{q}{q-1})^m \cdot \prod_{i=1}^m q^{-n+i\ell} = a \cdot (\tfrac{q}{q-1})^m \cdot q^{-mn+\frac{1}{2}\ell m(m+1)}.
\]
Now $\frac{n}{\ell}-1 < m \le \frac{n}{\ell}$. We conclude that
\[
P_G(w) \le a \cdot (\tfrac{q}{q-1})^{\frac{n}{\ell}} \cdot q^{-n(\frac{n}{\ell}-1) +\frac{1}{2}n(\frac{n}{\ell}+1)} = a \cdot (\tfrac{q}{q-1})^{\frac{n}{\ell}}
\cdot  q^{-\frac{n^2}{2\ell}+\frac{3n}{2}}.
\]
For $\ell \ge 3$, this is less than $q^{-\frac{n^2}{(2+\e)\ell}}$, hence less than $|G|^{-\frac{1}{(2+\e)\ell}}$, provided
$n \ge c(\e)\ell$, where $c(\e) = 4\,(1+\frac{2}{\e})$. And for $\ell = 2$, the same assertion holds using the upper bound for $i_2(G)$, the number of involutions in $G$, given by \cite[1.3]{LLS} (noting that for the word $w=x^2$, $P_G(w) = \frac{i_2(G)+1}{|G|}$).

This completes the proof of Theorem \ref{prob} for $G=GL_n(q)$, and the same argument replacing $G$ by $SL_n(q)$ gives the result for $SL_n(q)$.

Now let $G = Cl_n(q)$ be a classical group with natural module $V = (\F_{Q})^n$, where $Q=q^2$ if $G$ is unitary and $Q=q$ otherwise. Let $(\,,\,)$ be the associated bilinear or sesquilinear form on $V$ preserved by $G$, and when $G$ is orthogonal, let $R$ be the associated quadratic form. Assume $n > 4\ell$.

The proof is rather similar to the previous proof for $GL_n$.
Let $a_1,\ldots,a_k$ be free generators for $F_k$, and let $w = w_{\ell}\cdots w_1$, where each $w_i \in \{a_1^{\pm1},\ldots ,a_k^{\pm1}\}$. Let $g_1\ldots,g_k$ be a random sequence of elements of $G$.
Let $v_1 \in V$ be a nonzero singular vector, and define its trajectory $v_1^0,\ldots v_1^{\ell}$ as before.
Assume that $v_1^0,\ldots,v_1^{j-1}$ are linearly independent. Again, (\ref{rest}) holds, excluding at most $q^{j-1}$ possibilities for $v_1^j$. Moreover the values of $(v_1^j,v_1^i)$ are specified for the vectors $v_1^i$ for which $w_i=w_j$ or $w_{i+1}=w_j^{-1}$. Hence there are at least $Q^{n-j}-Q^{j-1}$ possibilities for $v_1^j$, and so
\[
{\bf P}\left(v_1^j \in {\rm Sp}(v_1^0,\ldots ,v_1^{j-1})\,|\,v_1^0,\ldots ,v_1^{j-1}\right) \le \frac{Q^j}
{Q^{n-j}-Q^{j-1}}.
\]
It follows that
\[
{\bf P}\left(v_1^0,\ldots ,v_1^\ell \hbox{ lin. dep.}\right) \le \sum_{j=1}^\ell \frac{Q^j}{Q^{n-j}-Q^{j-1}} \le
\frac{Q}{Q-1} \cdot \frac{Q^\ell}{Q^{n-\ell}-Q^{\ell-1}}.
\]
Now as before define further trajectories $v_i^0,\ldots, v_i^\ell$ for $1\le i \le m$, where $m<\frac{n}{2\ell}$. Arguing as above we obtain
\[
{\bf P}\left(v_i^0,\ldots ,v_i^{\ell} \hbox{ lin. dep.}\,|\, v_j^0,\ldots ,v_j^\ell \hbox{ for } j<i \right) \le
\frac{Q}{Q-1} \cdot \frac{Q^{i\ell}}{Q^{n-i\ell}-Q^{i\ell-1}}.
\]
If $w(g_1,\ldots,g_k)=1$, then $v_i^\ell = v_i$ for all $i$, and hence
\[
\begin{array}{ll}
P_G(w) & \le \prod_{i=1}^m {\bf P}\left(v_i^\ell = v_i\,|\,v_j^\ell=v_j \hbox{ for } 1\le j\le i-1\right) \\
  &  \le (\tfrac{Q}{Q-1})^m \prod_1^m \frac{Q^{i\ell}}{Q^{n-i\ell}-Q^{i\ell-1}} = (\tfrac{Q}{Q-1})^m \prod_1^m \frac{1}{Q^{n-2i\ell}-Q^{-1}}.
\end{array}
\]
Set $m=\lfloor \frac{n}{2\ell} \rfloor$. Arguing as above, this leads to
\[
P_G(w) \le a\cdot (\tfrac{Q}{Q-1})^{\frac{n}{2\ell}} Q^{-\frac{n^2}{4\ell}+\frac{3n}{2}}.
\]
As before, this gives $P_G(w) \le  |G|^{-\frac{1}{(2+\e)\ell}}$ provided $n \ge 7\,(1+\frac{2}{\e})\ell $.

This completes the proof of Theorem \ref{prob}.

\section{Deduction of Theorem \ref{girth}}

Let $G$ be a finite group and $S$ a sequence of $k$ random elements of $G$ chosen independently. For $\ell \ge 1$, define $P_G(\ell)$ to be the maximum of $P_G(w)$ over all words $w \in F_k$ of length $|w| = \ell$. Then as in \cite[Sec. 2]{GHSSV} by the well-known union bound, for any positive integer $L$ we have
\begin{equation}\label{eqn1}
{\bf P} (\hbox{girth}(\G(G,S)) \le L) \le  \sum_{|w|\le L} P_G(w) = \sum_{\ell=1}^{L} 2k(2k-1)^{\ell-1} P_G(\ell).
\end{equation}
Now let  $G = Cl_n(q)$ and choose $\e$ with $0<\e \le \log_{2k-1}q$. Then Theorem \ref{prob} gives $P_G(l) \le |G|^{-\frac{1}{(2+\e)\ell}}$ for $\ell \le \frac{n}{c}$, where $c =c(\e) = 7\,(1+\frac{2}{\e})$. Hence,
taking  $L \le \frac{n}{c}$,  the right hand side in (\ref{eqn1}) is bounded above by
\[
E:=\frac{k}{k-1} (2k-1)^{n/c} |G|^{-\frac{c}{(2+\e)n}}.
\]
Since $\frac{c}{2+\e} = \frac{7}{\e}$, we have
\[
\log_{2k-1}E \le 1+\frac{n}{c} - \frac{3.5(n-1)}{\e}\log_{2k-1}q,
\]
and by the choice of $\e$  this tends to 0 as $|G|\to \infty$. Hence the girth is at least $\frac{n}{c}$. Fixing $k$ and $\e$, this is of the order of $b(q)\sqrt{\log |G|}$, where $b(q) = (\log q)^{-\frac{1}{2}}$. Theorem \ref{girth} follows.

\section{Proof of Proposition \ref{DGbd}}

We first prove part (i). Fix $k \ge 2$. Let $G = G(q)$ be a simple group of Lie type of fixed rank and
let $S$ be a sequence of $k$ independently chosen random elements of $G$.  By \cite{LiSh1}, $S$ generates $G$
almost surely. By \cite[Thm. 4]{GHSSV}, the girth $g$ of $\G(G,S)$ satisfies
\begin{equation}\label{gbd}
g \ge c_1 \log |G|
\end{equation}
almost surely for some positive absolute constant $c_1$.
Let $T$ be the symmetric set consisting of the elements of $S$ and their inverses.
Write $h = \lfloor \frac{g-1}{2}\rfloor$. Since $g$ is the girth, we have $|T^h| \ge (2k-1)^h$ almost surely.
Let $A=T^h$. Then it follows from (\ref{gbd}) that $|A| \ge |G|^\d$ for some positive absolute constant $\d$ almost surely.

By the Product Theorem \cite{BGT, PS}, there is a positive absolute constant $\e$ such that for any symmetric generating subset $B$ of $G$,
either $B^3 = G$ or $|B^3| \ge |B|^{1+\e}$. It follows inductively that if $m$ is chosen minimally such that $\d (1+\e)^m\ge 1$, then we have
\[
A^{3^m} = G.
\]
Note that $m$ is an absolute constant. It follows that
\[
d = \hbox{diam}(\G(G,S)) \le  3^m \cdot h \le c_2 \cdot g
\]
almost surely, where $c_2 = 3^m/2$, as required.


Now we prove part (ii) of  Proposition \ref{DGbd}.
Let $G = Alt_n$ and let $S$ be a sequence of $k$ independently chosen random elements of $G$. It follows from \cite[Thm. 1.1]{Eber} that $g = \hbox{girth} (\G(G,S)) > c_1n^{1/3}$ almost surely, where $c_1$ is a positive absolute constant.  Also by \cite[Thm. 1.1]{HSZ}, $d = \hbox{diam}(\G(G,S)) < c_2 n^2(\log n)^{c_3}$. The conclusion follows.

The proof of part (iii) relies on \cite[Thm. 1.4]{BY}, showing that, for $G = Cl_n(q)$, the diameter $d$ of any
connected Cayley graph of $G$ satisfies
\[
d \le q^{O(n (\log n + \log q)^3)}.
\]
Combining this with the fact that the girth $g$ of $\G(G,S)$ satisfies $g \ge B(k)n$ almost surely (see the remark following Theorem \ref{girth}) we easily derive part (iii).

\section{Proof of Results \ref{pos}, \ref{posprop}, \ref{rep} and \ref{growth}}

The proof of \cite[Thm. 3]{GHSSV} given on p.106 was shown to contain an error by Eberhard \cite{Eber}. However, the error pointed out in \cite[Sec. 3]{Eber} only  pertains if there is a value of $i$ such that both $a_i$ and $a_i^{-1}$  occur in the word $w$. Hence, if we restrict to positive words, the inequality displayed as (6) on p.106 of \cite{GHSSV} holds. Proposition \ref{posprop} is just this bound.  Now Theorem \ref{pos} follows, just as in \cite[p.106]{GHSSV}.

To prove Corollaries \ref{rep} and \ref{growth} we may assume that $\Gamma = \la a_1, \ldots , a_k : w(a_1, \cdots , a_k) = 1 \ra$
where $w$ is the relation (resp. the positive relation) of minimal length $\ell$ (since our group is a quotient of the group above).

Note that, for an algebraic group $G$, the variety ${\rm Hom}(\Gamma,G)$ can be identified with the subvariety of $G^k$
defined by the equation $w(g_1, \ldots , g_k) = 1$.
The proof of Corollary \ref{rep} now follows using Theorem \ref{prob} and Lang-Weil estimates, as in \cite[Sec. 7]{LiSh3} and \cite[Sec. 4]{LS}.

Finally, to prove Corollary \ref{growth} we combine Proposition \ref{posprop} with the well-known inequality
\[
a_n(\Gamma) \le |{\rm Hom}(\Gamma, S_n)|/(n-1)! = P_{S_n}(w) n!^{k-1} \cdot n,
\]
which follows from \cite[1.1]{LuSe}.

\end{document}